
%
%
%
%
%
%
\magnification=\magstephalf      
%
%
\vsize=7.5truein                 
\hsize=5.2truein                 
\newskip\stdskip                 
\stdskip=6pt plus3pt minus3pt    
\medskipamount=\stdskip          
\parindent=0pt                   
\parskip=\stdskip                
\abovedisplayskip=\stdskip       
\belowdisplayskip=\stdskip       
\mathsurround=0.75pt             
\overfullrule=0pt                
%
%
\def\ppar{\par\goodbreak\vskip 8pt plus 4pt minus 4pt}     
%
%
\def\stdspace{\hskip 0.75em plus 0.15em\ignorespaces}
\let\qua\stdspace 
%
%
%
%
%
%
%
\def\hexnumber#1{\ifcase#1 0\or 1\or 2\or 3\or 4\or 5\or 6\or 7\or 8\or
 9\or A\or B\or C\or D\or E\or F\fi}
%
%
\font\thirtnmsa=msam10 scaled 1315    
\font\tenmsa=msam10          \font\ninemsa=msam9
\font\sevenmsa=msam7         \font\sixmsa=msam6
\font\fivemsa=msam5
%
%
\newfam\msafam                  \textfont\msafam=\tenmsa
\scriptfont\msafam=\sevenmsa    \scriptscriptfont\msafam=\fivemsa
\edef\hexa{\hexnumber\msafam}        
\def\msa{\fam\msafam\tenmsa}         
%
%
\font\thirtnmsb=msbm10 scaled 1315   
\font\tenmsb=msbm10      \font\ninemsb=msbm9
\font\sevenmsb=msbm7     \font\sixmsb=msbm6
\font\fivemsb=msbm5
%
\newfam\msbfam                   \textfont\msbfam=\tenmsb       
\scriptfont\msbfam=\sevenmsb     \scriptscriptfont\msbfam=\fivemsb
\edef\hexb{\hexnumber\msbfam}    
\def\msb{\fam\msbfam\tenmsb}     
%
%
\font\thirtneufm=eufm10 scaled 1315   
\font\teneufm=eufm10                 \font\nineeufm=eufm9
\font\seveneufm=eufm7                \font\sixeufm=eufm6
\font\fiveeufm=eufm5
%
\newfam\eufmfam                    \textfont\eufmfam=\teneufm
\scriptfont\eufmfam=\seveneufm     \scriptscriptfont\eufmfam=\fiveeufm
\edef\hexf{\hexnumber\eufmfam}      
\def\frak{\fam\eufmfam\teneufm}     
%
%
%
\font\thirtnrm=cmr10 scaled 1315    
\font\ninerm=cmr9                   \font\sixrm=cmr6   
%
\font\thirtni=cmmi10 scaled 1315    
\font\ninei=cmmi9                   \font\sixi=cmmi6  
%
\font\thirtnsy=cmsy10 scaled 1315   
\font\ninesy=cmsy9                  \font\sixsy=cmsy6  
%
\font\thirtnbf=cmbx10 scaled 1315   
\font\ninebf=cmbx9                  \font\sixbf=cmbx6  
%
%
\font\thirtnex=cmex10 scaled 1315   
\font\nineex=cmex9                  
%
%
\font\thirtnit=cmti10 scaled 1315  
\font\nineit=cmti9                  
%
\font\thirtnsl=cmsl10 scaled 1315  
\font\ninesl=cmsl9                  
%
\font\thirtntt=cmtt10 scaled 1315  
\font\ninett=cmtt9                  
%
%
%
%
\def\small{%
%
%
\textfont0=\ninerm \scriptfont0=\sixrm \scriptscriptfont0=\fiverm
\def\rm{\fam0\ninerm}
%
%
\textfont1=\ninei \scriptfont1=\sixi \scriptscriptfont1=\fivei
%
%
\textfont2=\ninesy \scriptfont2=\sixsy \scriptscriptfont2=\fivesy
%
%
\textfont3=\nineex \scriptfont3=\nineex \scriptscriptfont3=\nineex
%
%
\textfont\bffam=\ninebf \scriptfont\bffam=\sixbf
\scriptscriptfont\bffam=\fivebf \def\bf{\fam\bffam\ninebf}%
%
%
\textfont\itfam=\nineit \def\it{\fam\itfam\nineit}%
\textfont\slfam=\ninesl \def\sl{\fam\slfam\ninesl}%
\textfont\ttfam=\ninett \def\tt{\fam\ttfam\ninett}%
%
%
%
\textfont\msafam=\ninemsa \scriptfont\msafam=\sixmsa
\scriptscriptfont\msafam=\fivemsa \def\msa{\fam\msafam\ninemsa}%
%
%
\textfont\msbfam=\ninemsb \scriptfont\msbfam=\sixmsb
\scriptscriptfont\msbfam=\fivemsb \def\msb{\fam\msbfam\ninemsb}%
%
%
\textfont\eufmfam=\nineeufm  \scriptfont\eufmfam=\sixeufm
\scriptscriptfont\eufmfam=\fiveeufm \def\frak{\fam\eufmfam\nineeufm}%
%
%
%
\normalbaselineskip=11pt%
\setbox\strutbox=\hbox{\vrule height8pt depth3pt width0pt}%
%
%
\normalbaselines\rm
%
%
\stdskip=4pt plus2pt minus2pt    
\medskipamount=\stdskip          
\parskip=\stdskip                
\abovedisplayskip=\stdskip       
\belowdisplayskip=\stdskip       
\def\ppar{\par\goodbreak\vskip 6pt plus 3pt minus 3pt}%
%
%
\def\section##1{\global\advance\sectionnumber by 1
\vskip-\lastskip\penalty-800\vskip 20pt plus10pt minus5pt 
\egroup{\bf\number\sectionnumber\quad##1}\bgroup\small         
\vskip 6pt plus3pt minus3pt
\nobreak\resultnumber=1}
}    
%
\def\beginsmall{\bgroup\small}
\let\endsmall\egroup
%
%
%
%
\def\large{%
\textfont0=\thirtnrm \scriptfont0=\ninerm \scriptscriptfont0=\sevenrm
\def\rm{\fam0\thirtnrm}%
\textfont1=\thirtni \scriptfont1=\ninei \scriptscriptfont1=\seveni
\textfont2=\thirtnsy \scriptfont2=\ninesy \scriptscriptfont2=\sevensy
\textfont3=\thirtnex \scriptfont3=\thirtnex \scriptscriptfont3=\thirtnex
\textfont\bffam=\thirtnbf \scriptfont\bffam=\ninebf
\scriptscriptfont\bffam=\sevenbf \def\bf{\fam\bffam\thirtnbf}%
\textfont\itfam=\thirtnit \def\it{\fam\itfam\thirtnit}%
\textfont\slfam=\thirtnsl \def\sl{\fam\slfam\thirtnsl}%
\textfont\ttfam=\thirtntt \def\tt{\fam\ttfam\thirtntt}%
\textfont\msafam=\thirtnmsa \scriptfont\msafam=\ninemsa
\scriptscriptfont\msafam=\sevenmsa \def\msa{\fam\msafam\thirtnmsa}%
\textfont\msbfam=\thirtnmsb \scriptfont\msbfam=\ninemsb
\scriptscriptfont\msbfam=\sevenmsb \def\msb{\fam\msbfam\thirtnmsb}%
\textfont\eufmfam=\thirtneufm  \scriptfont\eufmfam=\nineeufm
\scriptscriptfont\eufmfam=\seveneufm \def\frak{\fam\eufmfam\teneufm}%
\normalbaselineskip=16pt%
\setbox\strutbox=\hbox{\vrule height11.5pt depth4.5pt width0pt}%
\normalbaselines\rm}%
\let\Large\large   
%
\def\Bbb#1{{\msb#1}}

%

\def\re{\Bbb R}
%
\mathchardef\plussquare="0\hexa01
\mathchardef\nge="3\hexb0B
\mathchardef\maltesecross="0\hexa7A
\mathchardef\del="0\hexf01
%
%
%
%
\font\sc=cmcsc10
%
%
%
%
\def\sqr#1#2{{\vcenter{\vbox{\hrule  height.#2truept
	\hbox{\vrule width.#2truept height#1truept 
	\kern#1truept \vrule width.#2truept}
	\hrule height.#2truept}}}}
\def\sq{\sqr55}    
%
%
%
%
\newcount\sectionnumber            
\newcount\resultnumber             
\sectionnumber=0\resultnumber=1    
%
%
%
\def\section#1{\global\advance\sectionnumber by 1
\xdef\nextkey{\number\sectionnumber}
\vskip-\lastskip\penalty-800\vskip 20pt plus10pt minus5pt 
{\large\bf\number\sectionnumber\quad#1}         
\vskip 8pt plus4pt minus4pt
\nobreak\resultnumber=1}      
%
%
%
%
%
\def\sh#1{\vskip-\lastskip\ppar{\bf #1}\par\nobreak\medskip}         
%
%
%
%

%
\def\proc#1{\xdef\nextkey{\number\sectionnumber.\number\resultnumber}%
\vskip-\lastskip\ppar\bf%
\noindent#1\ \number\sectionnumber.\number\resultnumber
\stdspace\sl\global\advance\resultnumber by 1\ignorespaces}
\def\endproc{\rm\ppar} 
%
%
\def\prf{\vskip-\lastskip\ppar\noindent{\bf Proof}%
\stdspace\rm}                            
\def\qed{\hfill$\sq$\par\goodbreak\rm}   
\def\endprf{\unskip\stdspace\hbox{}
\hfill$\sq$\par\medskip}                 
%
%
%
%
%
%
%
%
\def\proclaim#1{\vskip-\lastskip\ppar\bf%
\noindent#1\stdspace\sl\ignorespaces} 

%
%
%
%
\def\rk#1{\vskip-\lastskip\ppar{\bf #1}\stdspace\ignorespaces}                
\def\endrk{\par\medskip}
%
%
%
%
%
%
\def\label{\xdef\nextkey{\number\sectionnumber.\number\resultnumber}%
\number\sectionnumber.\number\resultnumber
\global\advance\resultnumber by 1}
%
%
%
%
%
%
%
%
%
%
%
%
%
%
%
%
\newcount\refnumber              
\refnumber=1                     
\long\def\reflist#1\endreflist{%
\long\def\thereflist{#1}{\def\refkey##1##2\par{\xdef##1{\number\refnumber}%
\global\advance\refnumber by 1}%
\def\key##1##2\par{\expandafter\xdef%
\csname##1\endcsname{\number\refnumber}%
\global\advance\refnumber by 1}#1\par}}
\long\def\references{%
\penalty-800\vskip-\lastskip\vskip 15pt plus10pt minus5pt 
{\large\bf References}\ppar 
{\leftskip=25pt\frenchspacing    
\small\parskip=3pt plus2pt       
\def\refkey##1##2\par{\noindent  
\llap{[##1]\stdspace}\ignorespaces##2\par}         
\def\key##1##2\par{\noindent  
\llap{[\ref{##1}]\stdspace}\ignorespaces##2\par}  
\def\,{\thinspace}\thereflist\par}}
%
%
%
\newcount\footnotenumber         
\footnotenumber=1                
\def\fnote#1{\xdef\nextkey{\number\footnotenumber}%
{\small\ifnum\footnotenumber>9\parindent=14pt%
\else\parindent=10pt\fi\footnote{$^{\number\footnotenumber}$}%
{\hglue-5pt#1}\global\advance\footnotenumber by 1}}
%
%
%
%
%
%
%
\newcount\figurenumber          
\figurenumber=1                 
\def\caption#1{\xdef\nextkey{\number\figurenumber}%
\cl{\small Figure \number\figurenumber: #1}%
\global\advance\figurenumber by 1}
\def\figurelabel{\xdef\nextkey{\number\figurenumber}%
\cl{\small Figure \number\figurenumber}%
\global\advance\figurenumber by 1}
\long\def\figure#1\endfigure{{\xdef\nextkey{\number\figurenumber}%
\let\captiontext\relax\def\caption##1{\xdef\captiontext{##1}}%
\midinsert\cl{\ignorespaces#1\unskip\unskip\unskip\unskip}\vglue6pt\cl{\small 
Figure \number\figurenumber\ifx\captiontext\relax\else: \captiontext
\fi}\endinsert\global\advance\figurenumber by 1}}
%
%
%
%
%
%
%
\def\nextkey{??}   
%
\def\key#1{\expandafter\xdef\csname #1\endcsname{\nextkey}}
\def\ref#1{\expandafter\ifx\csname #1\endcsname\relax
\immediate\write16{Reference {#1} undefined}??\else
\csname #1\endcsname\fi}
%
%
%
%
%
%
%
\newread\gtinfile
\newwrite\gtreffile
\def\useforwardrefs{
\openin\gtinfile\jobname.ref
\ifeof\gtinfile
\closein\gtinfile
\immediate\write16{No file \jobname.ref}
\else
\closein\gtinfile
\input \jobname.ref
\fi
\immediate\openout\gtreffile \jobname.ref
%
%
\def\key##1{{\def\\{\noexpand}%
\expandafter\xdef\csname ##1\endcsname{\nextkey}%
\immediate\write\gtreffile{\\\expandafter\\\def\\\csname ##1\\\endcsname%
{\nextkey}}}}
%
%
\long\def\reflist##1\endreflist{%
\long\def\thereflist{##1}{\def\refkey####1####2\par{\xdef####1{%
\number\refnumber}{\def\\{\noexpand}\immediate\write\gtreffile
{\\\def\\####1{\number\refnumber}}}\global\advance\refnumber by 1}%
\def\key####1####2\par{\expandafter\xdef%
\csname####1\endcsname{\number\refnumber}%
{\def\\{\noexpand}\immediate\write\gtreffile
{\\\expandafter\\\def\\\csname ####1\\\endcsname{\number\refnumber}}}
\global\advance\refnumber by 1}##1\par}}
\long\def\biblio##1\endbiblio{\reflist##1\endreflist\references}%
%
%
\def\numkey##1{{\def\\{\noexpand}%
\xdef##1{\number\sectionnumber.\number\resultnumber}
\immediate\write\gtreffile{\\\def\\##1%
{\number\sectionnumber.\number\resultnumber}}}}
\def\seckey##1{{\def\\{\noexpand}\xdef##1{\number\sectionnumber}
\immediate\write\gtreffile{\\\def\\##1{\number\sectionnumber}}}}
\def\figkey##1{\xdef##1{\number\figurenumber}%
{\def\\{\noexpand}\immediate\write\gtreffile%
{\\\def\\##1{\number\figurenumber}}}
\number\figurenumber\global\advance\figurenumber by 1}
}   
%
%
%
%
\def\figkey#1{\xdef#1{\number\figurenumber}%
\number\figurenumber\global\advance\figurenumber by 1}
\def\fig#1#2\endfig{%
\midinsert\cl{#2}\vglue6pt\cl{\small Figure #1}\endinsert}
\def\newfig{\number\figurenumber\global\advance\figurenumber by 1}
\def\numkey#1{\xdef#1{\number\sectionnumber.\number\resultnumber}}
\def\seckey#1{\xdef#1{\number\sectionnumber}}
%
%
%
%
%
%
%
%
%
\def\verb{\catcode`\"=\active}       
\def\brev{\catcode`\"=12}            
\brev                                
\verb                                
{\obeyspaces\gdef {\ }}              
{\catcode`\`=\active\gdef`{\relax\lq}}
\def"{%
\begingroup\baselineskip=12pt\def\par{\leavevmode\endgraf}%
\tt\obeylines\obeyspaces\parskip=0pt\parindent=0pt%
\catcode`\$=12\catcode`\&=12\catcode`\^=12\catcode`\#=12%
\catcode`\_=12\catcode`\~=12%
\catcode`\{=12\catcode`\}=12\catcode`\%=12\catcode`\\=12%
\catcode`\`=\active\let"\endgroup}
\brev      
%
%
%
%
%
%
\def\item#1{\par\leavevmode\llap{#1\stdspace}%
\ignorespaces}                             
%
%

%
%
\def\co{\colon\thinspace}    
\def\np{\vfil\eject}         
\def\nl{\hfil\break}         
\def\cl{\centerline}         
\def\agt{{\mathsurround=0pt\it$\cal A\mskip-.7mu$lgebraic \&\ 
$\cal G\mskip-2mu$eometric $\cal T\!\!$opology}}  
%
%
%

%
%
%
%
%
\def\title#1{\def\thetitle{#1}}

\def\author#1{\edef\previousauthors{\theauthors}
 \ifx\theauthors\relax\def\theauthors{#1}\else
 \def\theauthors{\previousauthors\par#1}\fi}

%
\def\address#1{\edef\previousaddresses{\theaddress}
 \ifx\theaddress\relax\def\theaddress{#1}\else
 \def\theaddress{\previousaddresses\par\vskip 2pt\par#1}\fi}
\def\secondaddress#1{\edef\previousaddresses{\theaddress}
 \ifx\theaddress\relax\def\theaddress{#1}\else
 \def\theaddress{\previousaddresses\par{\rm and}\par#1}\fi}   

\def\email#1{\edef\previousemails{\theemail}
 \ifx\theemail\relax\def\theemail{#1}\else
 \def\theemail{\previousemails\hskip 0.75em\relax#1}\fi}
\def\secondemail#1{\edef\previousemails{\theemail}
 \ifx\theemail\relax\def\theemail{#1}\else
 \def\theemail{\previousemails\hskip 0.75em{\rm and}\hskip 0.75em
 \relax#1}\fi}
\def\url#1{\edef\previousurls{\theurl}
 \ifx\theurl\relax\def\theurl{#1}\else
 \def\theurl{\previousurls\hskip 0.75em\relax#1}\fi}
\def\secondurl#1{\edef\previousurls{\theurl}
 \ifx\theurl\relax\def\theurl{#1}\else
 \def\theurl{\previousurls\hskip 0.75em{\rm and}\hskip 0.75em
 \relax#1}\fi}
\long\def\abstract#1\endabstract{\long\def\theabstract{#1}}
\def\primaryclass#1{\def\theprimaryclass{#1}}
\def\secondaryclass#1{\def\thesecondaryclass{#1}}
\def\keywords#1{\def\thekeywords{#1}}
%
%
\let\\\par\let\thetitle\relax\let\theshorttitle\relax
\let\theauthors\relax\let\theshortauthors\relax
\let\theaddress\relax\let\theshortaddress\relax
\let\theemail\relax\let\theurl\relax
\let\theabstract\relax\let\theprimaryclass\relax
\let\thesecondaryclass\relax\let\thekeywords\relax
%
%
%
%
\long\def\maketitlepage{    

\vglue 0.2truein   

%
{\parskip=0pt\leftskip 0pt plus 1fil\def\\{\par\smallskip}{\large
\bf\thetitle}\par\medskip}   

\vglue 0.15truein 

%
{\parskip=0pt\leftskip 0pt plus 1fil\def\\{\par}{\sc\theauthors}
\par\medskip}%
 
\vglue 0.1truein 

%
{\small\parskip=0pt
{\leftskip 0pt plus 1fil\def\\{\par}{\sl\theaddress}\par}
\ifx\theemail\relax\else  
\vglue 5pt \def\\{\stdspace{\rm and}\stdspace} 
\cl{Email:\stdspace\tt\theemail}\fi
\ifx\theurl\relax\else    
\vglue 5pt \def\\{\stdspace{\rm and}\stdspace} 
\cl{URL:\stdspace\tt\theurl}\fi\par}

\vglue 7pt 

{\bf Abstract}

\vglue 5pt

\theabstract

\vglue 7pt 

{\bf AMS Classification numbers}\quad Primary:\quad \theprimaryclass\par

Secondary:\quad \thesecondaryclass

\vglue 5pt 

{\bf Keywords:}\quad \thekeywords

\np  

}    
%
%
\long\def\makeshorttitle{    


%
{\parskip=0pt\leftskip 0pt plus 1fil\def\\{\par\smallskip}{\large
\bf\thetitle}\par\medskip}   

\vglue 0.05truein 

%
{\parskip=0pt\leftskip 0pt plus 1fil\def\\{\par}{\sc\theauthors}
\par\medskip}%
 
\vglue 0.03truein 

%
{\small\parskip=0pt
{\leftskip 0pt plus 1fil\def\\{\par}{\sl\ifx\theshortaddress\relax
\theaddress\else\theshortaddress\fi}\par}
\ifx\theemail\relax\else  
\vglue 5pt \def\\{\stdspace{\rm and}\stdspace} 
\cl{Email:\stdspace\tt\theemail}\fi
\ifx\theurl\relax\else    
\vglue 5pt \def\\{\stdspace{\rm and}\stdspace} 
\cl{URL:\stdspace\tt\theurl}\fi\par}

\vglue 10pt 


{\small\leftskip 25pt\rightskip 25pt{\bf Abstract}\stdspace\theabstract

{\bf AMS Classification}\stdspace\theprimaryclass
\ifx\thesecondaryclass\relax\else; \thesecondaryclass\fi\par
{\bf Keywords}\stdspace \thekeywords\par}
\vglue 7pt
}    
\let\maketitle\makeshorttitle        
%
%

\def\volumenumber#1{\def\thevolumenumber{#1}}
\def\volumeyear#1{\def\thevolumeyear{#1}}
\def\pagenumbers#1#2{\def\startpage{#1}\def\finishpage{#2}}
\def\published#1{\def\publishdate{#1}}
\def\received#1{\def\receiveddate{#1}}
\def\revised#1{\def\reviseddate{#1}}
\let\reviseddate\relax
\volumenumber{X}
\volumeyear{20XX}
\pagenumbers{1}{XXX}
\published{XX Xxxember 20XX}

\long\def\makeagttitle{   
\agt\hfill      
\hbox to 60truept{\vbox to 0pt{\vglue -14truept{\bf [Logo here]}\vss}\hss}
\break
{\small Volume \thevolumenumber\ (\thevolumeyear)
\startpage--\finishpage\nl
Published: \publishdate}

\vglue .2truein

{\parskip=0pt\leftskip 0pt plus 1fil\def\\{\par\smallskip}{\large
\bf\thetitle}\par\medskip}   
\vglue 0.05truein 

%
{\parskip=0pt\leftskip 0pt plus 1fil\def\\{\par}{\sc\theauthors}
\par\medskip}%
 
\vglue 0.03truein 


{\small\leftskip 25truept\rightskip 25truept{\bf Abstract}\stdspace\theabstract

{\bf AMS Classification}\stdspace\theprimaryclass
\ifx\thesecondaryclass\relax\else; \thesecondaryclass\fi\par
{\bf Keywords}\stdspace \thekeywords\par}\vglue 7truept

}   


\def\Addresses{\bigskip
{\small \parskip 0pt \leftskip 0pt \rightskip 0pt plus 1fil \def\\{\par}
\sl\theaddress\par\medskip \rm Email:\stdspace\tt\theemail\par
\ifx\theurl\relax\else\smallskip \rm URL:\stdspace\tt\theurl\par\fi}}

\def\agtart{
\hoffset 14truemm
\voffset 31truemm
\font\phead=cmsl9 scaled 950
\font\pnum=cmbx10 scaled 913
\font\pfoot=cmsl9 scaled 950
\headline{\vbox to 0pt{\vskip -4.5mm\line{\small\phead\ifnum
\count0=\startpage ISSN numbers are printed here
\hfill {\pnum\folio}\else\ifodd\count0\def\\{ }%
\ifx\theshorttitle\relax\thetitle\else\theshorttitle\fi\hfill{\pnum\folio}
\else\def\\{ and }{\pnum\folio}\hfill\ifx\theshortauthors\relax\theauthors
\else\theshortauthors\fi\fi\fi}\vss}}
\footline{\vbox to 0pt{\vglue 0mm\line{\small\pfoot\ifnum\count0=\startpage
Copyright declaration is printed here\hfill\else
\agt, Volume \thevolumenumber\ (\thevolumeyear)\hfill\fi}\vss}}
\let\maketitle\makeagttitle\let\makeshorttitle\makeagttitle}


\def\ifplaintex{\expandafter\ifx\csname documentclass\endcsname\relax}

\def\gtp{{\mathsurround=0pt\it $\cal G\mskip-2mu$eometry \&\ 
$\cal T\!\!$opology $\cal P\!$ublications}}  

\def\recd{{\small Received:\qua\receiveddate\ifx\reviseddate\relax
\else\qquad Revised:\qua\reviseddate\fi\par}} 


\def\lognumber#1{\def\thelognumber{#1}}
\def\volumenumber#1{\def\thevolumenumber{#1}}
\def\volumeyear#1{\def\thevolumeyear{#1}}
\def\papernumber#1{\def\thepapernumber{#1}}
\def\pagenumbers#1#2{\def\startpage{#1}\def\finishpage{#2}}
\def\published#1{\def\publishdate{#1}}

\def\received#1{\def\receiveddate{#1}}
\def\revised#1{\def\reviseddate{#1}}
\def\accepted#1{\def\accepteddate{#1}}

\def\asciiemail#1{\def\theasciiemail{#1}}

\long\def\asciiabstract#1{\long\def\theasciiabstract{#1}}
\def\asciikeywords#1{\def\theasciikeywords{#1}}


\let\\\par\let\thelognumber\relax\let\thevolumenumber\relax
\let\thepapernumber\relax\let\thevolumeyear\relax\let\startpage\relax
\let\finishpage\relax\let\publishdate\relax\let\receiveddate\relax
\let\reviseddate\relax\let\accepteddate\relax\let\theasciititle\relax
\let\theasciiauthors\relax
\let\theasciiabstract\relax\let\theasciikeywords\relax

\let\theasciiemail\relax


\ifplaintex
\font\logobig=cmssbx10 scaled 3836
\font\logomed=cmssbx10 scaled 2557
\else
\font\logobig=cmssbx10 scaled 4200
\font\logomed=cmssbx10 scaled 2800
\fi

\long\def\makeagttitle{   
\count0=\startpage
\agt\hfill      
\hbox to 45truept{\vbox to 0pt{\vglue -13truept{\logomed A\kern -.37em{\logobig 
T}\kern -.38em G}\vss}\hss}
\break
{\small Volume \thevolumenumber\ (\thevolumeyear)
\startpage--\finishpage\nl
Published: \publishdate}

\vglue .25truein

{\parskip=0pt\leftskip 0pt plus
1fil\def\\{\par\smallskip}{\Large\bf\thetitle}\par\medskip} \vglue
0.05truein

%
{\parskip=0pt\leftskip 0pt plus 1fil\def\\{\par}{\sc\theauthors}
\par\medskip}%
 
\vglue 0.03truein 


{\small\leftskip 25truept\rightskip 25truept{\bf Abstract}\stdspace\theabstract

{\bf AMS Classification}\stdspace\theprimaryclass
\ifx\thesecondaryclass\relax\else; \thesecondaryclass\fi\par
{\bf Keywords}\stdspace \thekeywords\par}\vglue 7truept

}   

\ifplaintex
\hoffset 14truemm
\voffset 31truemm
\font\phead=cmsl9 scaled 950
\font\pnum=cmbx10 scaled 913
\font\pfoot=cmsl9 scaled 950
\headline{\vbox to 0pt{\vskip -4.5mm\line{\small\phead\ifnum
\count0=\startpage ISSN 1472-2739 (on-line) 1472-2747 (printed)
\hfill {\pnum\folio}\else\ifodd\count0\def\\{ }%
\ifx\theshorttitle\relax\thetitle\else\theshorttitle\fi\hfill{\pnum\folio}
\else\def\\{ and }{\pnum\folio}\hfill\ifx\theshortauthors\relax\theauthors
\else\theshortauthors\fi\fi\fi}\vss}}
\footline{\vbox to 0pt{\vglue 0mm\line{\small\pfoot\ifnum\count0=\startpage
\copyright\ \gtp\hfill\else
\agt, Volume \thevolumenumber\ (\thevolumeyear)\hfill\fi}\vss}}
\else
\headsep 23pt
\footskip 35pt
\hoffset -4truemm
\voffset 12.5truemm
\font\lhead=cmsl9 scaled 1050
\font\lnum=cmbx10 
\font\lfoot=cmsl9 scaled 1050
\makeatletter
\def\@oddhead{{\small\lhead\ifnum\count0=\startpage ISSN 1472-2739 
(on-line) 1472-2747 (printed)\hfill {\lnum\number\count0}\else\ifodd\count0
\def\\{ }\ifx\theshorttitle\relax \thetitle \else\theshorttitle\fi\hfill
{\lnum\number\count0}\else\def\\{ and }{\lnum\number\count0}
\hfill\ifx\theshortauthors\relax 
\theauthors\else\theshortauthors\fi\fi\fi}}\def\@evenhead{\@oddhead}
\def\@oddfoot{\small\lfoot\ifnum\count0=\startpage\copyright\ \gtp\hfill\else
\agt, Volume \thevolumenumber\ (\thevolumeyear)\hfill\fi}
\def\@evenfoot{\@oddfoot}
\makeatother
\fi
\let\maketitlepage\makeagttitle
\let\makeshorttitle\maketitlepage
\let\maketitle\maketitlepage


\newwrite\gtoutfile
\long\gdef\makeheadfile{  
{\def\\{, }\def\s{ }
\immediate\openout\gtoutfile head.xxx
\immediate\write\gtoutfile{To: math@arxiv.org}
\immediate\write\gtoutfile{Subject: put OR rep NNNNN:ppppp}
\immediate\write\gtoutfile{--text follows this line--}
\immediate\write\gtoutfile{Proxy-for: \ifx\theasciiauthors\relax
\theauthors\else\theasciiauthors\fi\s<\ifx\theasciiemail\relax\theemail\else\theasciiemail\fi>}
\immediate\write\gtoutfile{\noexpand\\}
\immediate\write\gtoutfile{Authors: \ifx\theasciiauthors\relax
\theauthors\else\theasciiauthors\fi}
{\def\\{ }\immediate\write\gtoutfile{Title: \ifx\theasciititle\relax
\thetitle\else\theasciititle\fi}}
\immediate\write\gtoutfile{Subj-class: GT or SG, GR etc}
\immediate\write\gtoutfile{MSC-class: \theprimaryclass\ifx\thesecondaryclass\relax\else, \thesecondaryclass\fi}
\immediate\write\gtoutfile{Journal-ref: Algebr. Geom. Topol. \thevolumenumber\s
(\thevolumeyear) \startpage-\finishpage}
\immediate\write\gtoutfile{Comments: Published by Algebraic and
Geometric Topology at}
\immediate\write\gtoutfile{\s\s\s  http://www.maths.warwick.ac.uk/agt/AGTVol\thevolumenumber/agt-\thevolumenumber-\thepapernumber.abs.html}
\immediate\write\gtoutfile{\noexpand\\}
\immediate\write\gtoutfile{}
\ifx\theasciiabstract\relax
\immediate\write\gtoutfile{\theabstract}\else
\immediate\write\gtoutfile{\theasciiabstract}\fi
\immediate\write\gtoutfile{}
\immediate\write\gtoutfile{\noexpand\\}
\immediate\write\gtoutfile{}
\immediate\closeout\gtoutfile}}  

\def\maketitlepage{\makeagttitle\makeheadfile}
\let\makeshorttitle\maketitlepage
\let\maketitle\maketitlepage

\lognumber{28}
\volumenumber{3}
\volumeyear{2003}
\papernumber{28}
\published{25 September 2003}
\pagenumbers{857}{872}
\received{31 January 2003}
\revised{16 September 2003}
\accepted{24 September 2003}

\input xyall           

%
\def\squarediag#1#2#3#4#5#6#7#8{%
\xymatrix{
#1\ar[r]^-{#2}\ar[d]_{#4}&#3\ar[d]^{#5} \\
#6\ar[r]_-{#7}&#8}}
\def\upsquarediag#1#2#3#4#5#6#7#8{%
\xymatrix{
#1\ar[r]^-{#2}&#3 \\
#6\ar[r]_-{#7}\ar[u]^{#4}&#8\ar[u]_{#5}}}
\def\commdiag#1#2#3#4#5#6{%
$$\xymatrix{
#1\ar[r]^-{#2}\ar[dr]_{#4}&#3\ar[d]^{#5} \\
&#6}$$}
%
%
\def\tw{{\tt \char'176}}  
\def\E{{\cal E}}  
\def\I{{\cal I}}

\def\ep{\varepsilon}
\def\a{\longrightarrow}
\font\spec=cmtex10 scaled 1095 
\def\d{\hbox{\spec \char'017\kern 0.05em}} 
\def\inv{^{-1}}

\def\verts{\thinspace\vert\thinspace}

%
%
\reflist

\refkey\Boa
{\bf J\,M Boardman}, {\it Singularities of differentiable maps},
Publ. IHES, {33} (1967) 21--57

\refkey\Eliash
{\bf Y\,M Eliashberg}, {\it Surgery of singularities of smooth
mappings}, Izv. Akad. Nauk SSSR, Ser. Mat. Tom 36 (1972) no. 6,
English translation: Math. USSR-Izv. 6 (1972), 1302--1326

\refkey\EliGr
{\bf Y\,M Eliashberg}, {\bf M Gromov}, {\it Construction of
non-singular isoperimetric films}, Proc.  Steklov Math. Inst. 116
(1971) 13--28

\refkey\FRS
{\bf R Fenn}, {\bf C Rourke}, {\bf B Sanderson}, {\it James bundles
and applications},\nl {\tt http://www.maths.warwick.ac.uk/\tw
cpr/ftp/james.ps}

\refkey\Gr
{\bf M Gromov}, {\it Partial differential relations}, Springer--Verlag (1986)

\refkey\Hi
{\bf M Hirsch}, {\it Immersions of manifolds}, Trans. Amer. Math. Soc.
{93} (1959) 242--276

\refkey\Ja 
{\bf I James}, {\it Reduced product spaces}, Annals of Math. {62} (1955) 
170--197     

\refkey\Kos 
{\bf U Koschorke}, {\it Vector fields and other vector bundle morphisms --
a singularity approach}, Springer Lecture Notes Series, 847 (1981)

\refkey\KSb 
{\bf U Koschorke}, {\bf B Sanderson}, {\it Self-intersections and
higher Hopf invariants}, Topology, {17} (1978) 283-290

\refkey\LS 
{\bf R Lashof}, {\bf S.Smale}, {\it Self-intersections of immersed
manifolds}, J. Math. Mech. {8} (1959) 143-157

\refkey\Ma 
{\bf J\,P May}, {\it The geometry of iterated loop spaces},
Springer Lectures Notes Series, {271} Springer-Verlag (1972)

\refkey\Mi 
{\bf R\,J Milgram}, {\it Iterated loop spaces}, Annals of Math. 
{84} (1966) 386--403

\refkey\Config
{\bf Colin Rourke}, {\bf Brian Sanderson}, {\it Equivariant configuration
spaces}, J. London Math. Soc. 62 (2000) 544-552,
{\tt arXiv:math.GT/9711216}

\refkey\CompI
{\bf Colin Rourke}, {\bf Brian Sanderson}, {\it The compression
theorem I}, Geometry and Topology 5 (2001) 399--429, {\tt 
arXiv:math.GT/9712235}

\refkey\CompII
{\bf Colin Rourke}, {\bf Brian Sanderson}, {\it The compression theorem II:
directed embeddings}, Geometry and Topology 5 (2001) 431--440,
{\tt arXiv:math.GT/0003026}

\refkey\ImmThy
{\bf Colin Rourke}, {\bf Brian Sanderson}, {\it The compression
theorem}, pre-publication version {\tt arXiv:math.GT/9712235v2}

\refkey\Se 
{\bf G Segal}, {\it Configuration-spaces and iterated loop-spaces},
Invent. Math. {21} (1973) 213--221

\refkey\Sm 
{\bf S Smale}, {\it The classification of immersions of spheres in
Euclidean spaces}, Annals of Math. {69} (1959) 327--344

\refkey\Var
{\bf A Varley}, {\it Controlling tangencies from the global
viewpoint}, PhD thesis, Univ. of Warwick (2003)

\endreflist
%
%
\title{The compression theorem III: applications}

\author{Colin Rourke\\Brian Sanderson}  

\address{Mathematics Institute, University of Warwick\\
Coventry, CV4 7AL, UK}
\asciiemail{cpr@maths.warwick.ac.uk, bjs@maths.warwick.ac.uk}
\email{cpr@maths.warwick.ac.uk}
\secondemail{bjs@maths.warwick.ac.uk}
\url{http://www.maths.warwick.ac.uk/\char'176cpr\qua 
{\rm and}\qua\char'176bjs}

\abstract 
This is the third of three papers about the {\it Compression
Theorem\/}: if $M^m$ is embedded in $Q^q\times\re$ with a normal
vector field and if $q-m\ge1$, then the given vector field can be {\sl
straightened} (ie, made parallel to the given $\re$ direction) by an
isotopy of $M$ and normal field in $Q\times\re$.

The theorem can be deduced from Gromov's theorem on directed
embeddings [\Gr; 2.4.5 $(\rm C')$] and the first two parts gave
proofs.  Here we are concerned with applications.

We give short new (and constructive) proofs for immersion theory and
for the loops--suspension theorem of James et al and a new approach to
classifying embeddings of manifolds in codimension one or more, which
leads to theoretical solutions.

We also consider the general problem of controlling
the singularities of a smooth projection up to $C^0$--small isotopy
and give a theoretical solution in the codimension $\ge1$ case.
\endabstract

\asciiabstract{ 
This is the third of three papers about the Compression Theorem: if
M^m is embedded in Q^q X R with a normal vector field and if q-m > 0,
then the given vector field can be straightened (ie, made parallel to
the given R direction) by an isotopy of M and normal field in Q x R.
The theorem can be deduced from Gromov's theorem on directed
embeddings [Partial differential relations, Springer-Verlag (1986);
2.4.5 C'] and the first two parts gave proofs.  Here we are concerned
with applications.  We give short new (and constructive) proofs for
immersion theory and for the loops-suspension theorem of James et al
and a new approach to classifying embeddings of manifolds in
codimension one or more, which leads to theoretical solutions.  We
also consider the general problem of controlling the singularities of
a smooth projection up to C^0-small isotopy and give a theoretical
solution in the codimension >0 case.}

\primaryclass{57R25, 57R27, 57R40, 57R42, 57R52}

\secondaryclass{57R20,\break 57R45, 55P35, 55P40, 55P47}

\asciikeywords{Compression, embedding, isotopy, immersion, singularities, 
vector field, loops-suspension, knot, configuration space}

\maketitle

\section{Introduction}  
  
We work throughout in the smooth ($C^\infty$) category.  
The tangent bundle of a manifold $W$ is denoted $TW$ and the tangent
space at $x\in W$ is denoted $T_xW$. Throughout the paper,
``normal'' means independent (as in the usual meaning of ``normal
bundle'') and not necessarily perpendicular.

This is the third of a set of three papers about the following result:

\proclaim{Compression Theorem}Suppose that $M^m$ is embedded in 
$Q^q\times\re$ with a normal vector field and suppose that
$q-m\ge1$. Then the vector field can be {\sl straightened} (ie,
made parallel to the given $\re$ direction) by an isotopy of $M$ and
normal field in $Q\times\re$.\endproc

Thus the theorem moves $M$ to a position where it projects
by vertical projection (ie ``compresses'') to an immersion 
in $Q$.

The theorem can be deduced from Gromov's theorem on directed
embeddings [\Gr; 2.4.5 $(\rm C')$].  Proofs are given in parts I and
II [\CompI, \CompII].  Two immediate applications were given in
[\CompI; corollaries 1.1 and 1.2]; here we give more substantial
applications.

Immersion theory [\Hi, \Sm] implies the embedding is {\it regularly
homotopic} to an immersion which covers an immersion in $Q$ and using
configuration space models of multiple-loops-suspension spaces [\Ja,
\Ma, \Mi, \Se] it can be seen that the embedding is {\it bordant} (by a
bordism mapping to $M$) to an embedding which covers an immersion in
$Q$, see [\KSb].  Thus the new information which the compression
theorem provides is that the embedding is {\it isotopic} to an
embedding covering an immersion in $Q$.  Moreover we can apply the
compression theorem to give short and constructive proofs for both
immersion theory and configuration space theory.
  
The Compression Theorem can also be used to give a new approach to the
embeddings and knot problems for manifolds in codimension one or more,
which leads to theoretical solutions to both problems. 

The theorem also sheds light on the following problem:

\proclaim{$C^0$--Singularity Problem}Given $M\subset W$ and $p\,\co 
W\to Q$ a submersion, how much control do we have over the
singularities of $p|M$ if we are allowed a $C^0$--small isotopy of $M$
in $W$?\endproc

The problem includes the problem of controlling the singularities of a
map $f\co M\to Q$ by a $C^0$--small homotopy.  This is because we can
always factor $f$ as $f\times q\,\co M\subset Q\times\re^t\buildrel
{\rm proj}\over\a Q$ where $q\,\co M\to \re^t$ is an
embedding.

The Compression Theorem gives a necessary and sufficient condition for
desingularising the projection in the case that ${\rm dim}(W)={\rm
dim}(Q)+1$ and ${\rm dim}(M)<{\rm dim}(Q)$, namely that there should
exist an appropriate normal line field, and it can be extended to give
necessary and sufficient conditions for singularities of almost any
pre-specified type, namely that there should exist a line field with
those singularities.  Both these results extend to the case with just
the hypothesis ${\rm dim}(M)<{\rm dim}(Q)$ where ``line field'' is
replaced by ``plane field''.  They also have natural relative and
parametrised versions.

It is worth contrasting these results with the $C^\infty$--singularity
problem (``$C^0$--small'' is replaced by ``$C^\infty$--small'') where
there is the classical Thom--Boardman classification of
$C^\infty$--stable singularities [\Boa].  For most topological
purposes (eg for applications to homotopy theory) the $C^0$
classification is more natural than the classical $C^\infty$
classification.  Furthermore essential singularities (up to $C^0$
homotopy) have natural interpretations as generalised bordism
characteristic classes similar to those investigated by Korschorke
[\Kos] (see comments at the end of section 5).

The application to Immersion Theory is in section 2 and the
application to configuration space theory in section 3.  Section 4
contains the new approach to the embedding and knot problems and
section 5 is about the $C^0$--Singularity Theorem mentioned above.

\rk{Acknowledgement}We are extremely grateful to both Yasha Eliashberg 
and the referee for helpful comments on previous versions of this and
the other papers in this series.

\section{Immersion theory}

In the compression theorem, the existence of the immersion
of $M$ in $Q$ follows from immersion theory; however immersion
theory gives us no explicit information about this immersion, which is
only determined up to regular homotopy.  By contrast the
compression theorem gives us an explicit description of the
immersion in terms of the given embedding and normal vector field in
$Q\times\re$.  Moreover the compression theorem can be used
to give a new proof for immersion theory as we now show.  

Let $M$ be an $n$--manifold.  We shall explicitly describe a way of 
{\sl rotating the fibres} of 
the tangent bundle $TM$ into $M$.   Regard the zero section $M$ 
as `vertical' and the fibres as `horizontal'.   Consider
$TM$ as a smooth $2n$--manifold, then its tangent bundle 
restricted to $M$ is
the Whitney sum $TM\oplus TM$.  The two copies of $TM$ are 
the vertical copy parallel to $M$ and the horizontal copy
parallel to the fibres of $TM$.  Each vector $v\in TM$ then determines
two vectors $v_v$ and $v_h$ in $TM\oplus TM$ which span a plane.
In this plane we can `rotate' $v_h$ to $v_v$.  Since we are not
at this moment considering a particular metric `rotation' needs 
to be defined: to be precise we
consider the family of linear transformations of this plane given by 
$$v_h\mapsto cv_h+sv_v,\quad v_v\mapsto cv_v-sv_h\qquad \hbox{where}
\numkey\FRotate\leqno{\bf\label}$$ 
$$c=\cos{\textstyle {\pi\over2}}t,\quad s=\sin{\textstyle {\pi\over2}}t, 
\quad 0\le t\le1.$$  This
formula (applied to each such plane) determines a bundle
isotopy (a 1--parameter family of bundle isomorphisms) 
which is the required rotation of the fibres of $TM$ into $M$.  

\numkey\ImTh\proc{Immersion Theorem}Suppose that we 
are given a bundle monomorphism
$f\co TM\to TQ$ (ie, a map $M\to Q$ covered by a vector space 
monomorphism on each fibre) and that either $q-m\ge1$ or 
$q-m\ge0$ and each component of $M$ has relative boundary.
Then the restriction $f\vert\co M\to Q$ is homotopic to an 
immersion.\rm

\prf Composing $f$ with an exponential map for $TQ$
gives a map $g\co TM\to Q$ which embeds
fibres into $Q$.  Choose an embedding $q\co M\to \re^n$ for some $n$
and also denote by $q$ the map $TM\to \re^n$ given by projecting
the bundle $TM$ onto $M$ (the usual bundle projection) and then
composing with $q$.  We then have the embedding
$g\times q\co TM\to Q\times\re^n$.  The fibres of $TM$ are embedded
parallel to $Q$ and the $n$ directions parallel
to the axes of $\re^n$ determine $n$ independent vector fields
at $M$ normal to the fibres of $TM$.

Now choose a complement for $T(TM)\verts M\cong TM\oplus TM$ in
$T(Q\times\re)\verts M$.  Then the rotation of the fibres of $TM$ into
$M$ (formula \FRotate\ above) extends (by the identity on the
complement) to a bundle isotopy of $T(Q\times\re)\verts M$ which
carries the these $n$ fields normal to $M$ to yield $n$ independent
normal fields.  The result now follows from the multi-compression
theorem in part I [\CompI; corollary 4.5].\qed\ppar

The proof of the immersion theorem just given is very explicit, which
contrasts with the standard Hirsch--Smale approach [\Hi, \Sm] or the
proofs given by Gromov [\Gr].  Given a particular bundle monomorphism
$TM\to TQ$ the proof can be used to construct a homotopic immersion
$M\to Q$.  The only serious element of choice in the proof is the
embedding of $M$ in $\re^n$.  It is worth remarking that Eliashberg
and Gromov [\EliGr; Theorem 4.3.4] have also given a short proof of
immersion theory which yields an explicit immersion.

Notice that, in the proof just given, there was an explicit bundle
homotopy of the $n$ independent normal fields to $M$ to the vertical
fields.  Thus using the Normal Deformation Theorem [\CompI; Theorem
4.7] in place of the multi-compression theorem, gives a parametrised
version of the theorem.  It is easy to deduce the the usual statement
of immersion theory from this version.  

The parametrised proof can also be made very explicit by examining
carefully the above proof; for full details see [\ImmThy; section 6,
page 31].

\section{Loops--suspension theorem}

We next show that the compression theorem can be used to give a short
new proof of the classical result of James [\Ja] on the homotopy type
of loops--suspension and of the generalisation due to May [\Ma] and
Segal [\Se] and implicit in Milgram [\Mi].  In [\Config] the arguments
in this section are extended to both the equivariant case and the
disconnected case (where group completions are needed).

We denote the free topological monoid on a based space $X$ by
$X_\infty$ and denote the loop space on and suspension of $X$ by
$\Omega(X)$, $S(X)$, respectively.  We assume spaces are compactly
generated and Hausdorff and we assume base points non-degenerate. In
particular this means (up to homotopy) we can assume based spaces have
{\sl whiskers}, ie, the base point has a neighbourhood homeomorphic to
$\{1\}$ in the interval $[0,1]$. It will be convenient to identify the
whisker with $[0,1]$.

There is a map $k_X\co X_\infty\to\Omega S(X)$ 
defined in [\Ja].
Briefly, what $k_X$ does is to map the word $x_1\cdots x_m$
to a loop in $S(X)$ which comprises $m$ vertical loops passing
through $x_1, \dots, x_m$ respectively, with the time parameters
carefully adjusted to make the time spent on a subloop
go to zero as the corresponding point $x_i$ moves to the basepoint
of $X$. 

\numkey\LSone\proc{Theorem}Let $X$ be path-connected. 
Then $k_X\co X_\infty\to\Omega S(X)$
 is a weak homotopy equivalence.\rm

\prf There is a well-known equivalent definition of $X_\infty$ (up to
homotopy type) as the configuration space $C_1(X)$
of points in $\re^1$ labelled in $X$, of which we briefly recall the
definition.  Consider finite subsets of $\re^1$
labelled by points of $X$.  An equivalence relation on such
subsets is generated by deleting points labelled by the basepoint.
$C_1(X)$ is the set of equivalence classes.  The topology
is induced from the topologies of $\re^1$ and $X$.

We shall give a geometric description for the homomorphism 
$\pi_n(C_1X)\to \pi_n(\Omega SX)$ induced by $k_X$.  The construction 
is similar to a construction in [\KSb], which used transversality
to the base of a Thom space.  Here we use transversality to $X$ in
$X\times\re$ and to an interior point of the whisker in $X$.

We start by giving a geometric description for a map into 
$C_1(X)$.  Let $Q$ be a smooth manifold and $f\co Q\to C_1(X)$ a
map.  Then $f$ determines a subset of $Q\times\re^1$, continuously
labelled by $X$, such that at points labelled in $X-\{*\}$ the
projection on $Q$ is a local homeomorphism.  This subset is
determined up to introduction and deletion of points labelled
by $*$.  By approximating
the local map to $\re^1$ by a smooth map, we can assume, by a small
homotopy of $f$, that the subset of points labelled in $X-\{*\}$
is a smooth submanifold of $Q\times\re$.  Further, by using
transversality of the labelling map 
to an interior point $1\over2$ in the whisker in $X$, 
and composing with the stretch $[0,{1\over2}]\to[0,1]$, we may
assume that this labelled subset is in fact a smooth submanifold
$W$ with boundary, such that the boundary is the subset labelled by
$*$ and such that the projection on $Q$ is a local embedding
(of codimension 0).   Conversely, such a subset determines a
map $f\co Q\to C_1(X)$.  Notice that $W$ is canonically
framed in $Q\times\re^1$ by the $\re^1$--coordinate and also notice that 
if $Q$ and $f$ are based then $W$ can be assumed to be 
empty over the basepoint of $Q$. 

Next we interpret maps in $\Omega SX$.   A map $Q\to\Omega SX$
determines a map $f\co Q\times \re^1\to SX$.  By making $f$ transverse
to the suspension line through $1\over2$ in the whisker in $X$, 
and composing with the stretch
$[0,{1\over2}]\to[0,1]$, we may assume that $\overline{f\inv(SX-\{*\})}$
is a codimension 0 submanifold whose boundary maps to $*_{SX}$.  By
further making $f$ transverse to $X\times \{0\}$ we may assume that
$f\inv((X-\{*\})\times \{0\})$ is a framed codimension 1 submanifold 
$W$ with boundary, equipped with a map $l\co W\to X$ such that 
$l^{-1}*=\d W$ and that $f$ maps framing lines to suspension lines and the
rest is mapped to $*_{SX}$.  Conversely such a framed submanifold
determines a map $Q\to\Omega SX$.  Further if $Q$ and $f$ are based 
then $W$ is empty over the basepoint of $Q$.

With these geometric descriptions there is an obvious forgetful
map $[Q,C_1X]\to[Q,\Omega SX]$ which may be seen to be induced
by $k_X$.  Now consider the case when $Q=S^n$ and
consider a framed manifold $W$ representing an 
element of $\pi_n(SX)$. If it has a closed component we can 
change the labelling function to map a small disc to $*_X$
by using a collar on the disc and mapping the collar lines to
a path to the basepoint in $X$. Then the interior of the disc,
which is now labelled by $*_X$, may be deleted from $W$.
After eliminating all closed components in this way the 
compression theorem (the codimension 0 case) implies that
$\pi_n(C_1X)\to \pi_n(\Omega SX)$ is surjective.  

Injectivity is proved similarly by using the case $Q=S^n\times I$
and working relative to $Q\times\{0,1\}$.  \qed\ppar

Notice that only the Global Compression Theorem [\CompI; 2.1] (and addendum
(i)) were used for the proof of \LSone\ so the complete new proof is
short.  However, for the full loops--suspension theorem below we need
the multi-compression theorem.  (For the shortest proof of the full
theorem, use the short proof of the multi-compression theorem given in
[\CompII].)

Let $C_n(X)$ denote the configuration space of points
in $\re^n$ labelled in $X$, with, as before, points labelled
by $*$ removable.  There is a map $q_X\co C_n(X)\to
\Omega^n S^n(X)$ defined in a similar way to $k_X$ (as interpreted
in the last proof) using little
cubes around the points in $\re^n$ with axes parallel to the
axes of $\re^n$.

\proc{Theorem}Let $X$ be a path-connected topological space with a
non-degenerate basepoint then
$q_X$ is a weak homotopy equivalence.\rm

\prf  The proof is similar to the proof for theorem \LSone.
A map $Q\to C_n(X)$ can be see as partial cover of $Q$
embedded in $Q\times\re^n$ and a map $Q\to \Omega^n S^n(X)$ can
be see as a framed codimension $n$ submanifold of $Q\times\re^n$.

There is then a function Maps$(Q\to C_n(X))$ to Maps$(Q\to 
\Omega^n S^n(X))$ given by taking the parallel framing on the
partial cover.  This can be seen to be given by composition
with $q_X$.  Closed components in the submanifold of $Q\times\re^n$ 
may be punctured as in the last proof.  The codimension 0 case of
the multi-compression theorem now implies that $q_X$ 
induces a bijection between the sets of homotopy classes. \qed

\section{Embeddings and knots}

Two basic problems of differential topology are the {\it embedding
problem}: given two manifolds $M$ and $Q$, decide whether $M$ embeds
in $Q$, and the {\it knot problem}: classify embeddings of $M$ in $Q$
up to isotopy.  The corresponding problems for immersions (replace
embedding by immersion and isotopy by regular homotopy) are, in some
sense, solved by immersion theory; ie, solved by reducing
to vector bundle
problems for which there are standard obstruction theories.  Now if
$M$ embeds in $Q$ then it certainly immerses in $Q$ and if two
embeddings are isotopic then they are certainly regularly homotopic.
Thus it makes sense to consider the {\it relative embedding problem}:
decide whether a given immersion of $M$ in $Q$ is regularly homotopic
to an embedding, and the {\it relative knot problem}: given a regular
homotopy between embeddings, decide if it can be deformed to an
isotopy.  Since a manifold often immerses in a considerably lower
dimension than that in which it embeds, it makes sense to consider the
following more general problems.

\proc{Embedding problem}Suppose given an immersion $f\co
M\to Q$ and an integer $n\ge0$.  Decide whether $f\times0\,\co M\to 
Q\times\re^n$ is regularly homotopic to an embedding.\rm

\proc{Knot problem}Classify, up to isotopy, embeddings of $M$ in
$Q\times\re^n$ within the regular homotopy class of $f\times0\,\co M\to
Q\times\re^n$, where $f\co M\to Q$ is a given immersion.\rm
\ppar

The compression theorem gives substantial information on these
problems.  In particular we can give formal solutions which make it
easy to define obstructions to the existence of such embeddings or
isotopies.

\sh{The embedding theorem}

Turning first to the embedding problem, we have the following
application of the compression theorem.

\proc{Proposition}\key{CTProp}Suppose we are given an immersion
$f\co M\to Q$ and that $c=q-m\ge2$.  Suppose the immersion
$f\times0\,\co M\to Q\times\re^n$ is regularly homotopic to an
embedding $g$. Then $f\co M\to Q$ is regularly homotopic to an
immersion $f_1$ covered by an embedding isotopic to $g$ (ie an
embedding $g_1\co M\to Q\times\re^n$ such that $g_1(x)=(f_1(x),0)$ for
$x\in M$).

\prf 
 Let $F\co M\times I\to Q\times\re^n\times I$ be determined by the
regular homotopy of $f\times 0$ to $g$.  Thus $F_0=F\verts
(M\times\{0\}) = f\times0$ and $F_1=F\verts (M\times\{1\})=g$.  The
canonical $n$--frame on $F_0$ extends to an $n$--frame on $F$.  Apply
the multi-compression theorem to $F_1$ to yield an isotopy of $F_1$ to
$F_1'=g_1$ which compresses to an immersion $f_1\co M\to Q$.  Extend
the isotopy of $F_1$ to a regular homotopy of $F$ to $F'$ say, rel
$F_0$, and apply the parametrised multi-compression theorem [\CompI]
again to compress $F'$ to a regular homotopy between $f$ and
$g_1$.\break\hbox{}
\qed\ppar

Recall from section 3 that $C_n(X)$ denotes the configuration space of
points in $\re^n$ labelled in $X$, with points labelled by $*$
removable.

Suppose that we are given an immersion $f\co M\to Q$ covered by an
embedding $e\co M\to Q\times\re^n$.  Under these circumstances
Koschorke and Sanderson [\KSb] define a classifying map
$\alpha_f^n\co Q\to C_n(MO_c)$ where $MO_c$ is the Thom space of the
universal $c$--vector bundle $\gamma^c/BO_c$.  The definition is
roughly as follows.  The normal (disc) bundle $\xi$ on the immersion
$f$ immerses in $Q$ and each point $q\in Q$ then lies in the image of
a number of the fibres.  This defines a configuration in $\re^n$ 
by considering the covering immersion.  The labelling in $MO_c$ comes
from using the classifying map for $\xi$.

In detail the embedding $e$ has a normal bundle which splits a trivial
$n$--plane bundle say $\nu=\xi\oplus\epsilon^n$ and the embedding
extends to an embedding $e'\co D(\nu)\to Q\times\re^n$ so that
$D(\xi)\subset D(\nu){\buildrel proj\over\a} Q$ gives a tubular
neighbourhood of $M$ in $Q$.  We can assume that, for each $q\in Q$,
$\{q\}\times\re^n \cap e'(D(\xi))$ is finite.  Let $h\co D(\xi)\to
D(\gamma^c)$ be a classifying (disc) bundle map, where $\gamma^c/BO_c$
is a universal vector bundle.  Let $MO_c$ denote the Thom space
$T(\gamma^c)$.  Then we have a map $\alpha_f^n\co Q\to C_n(MO_c)$
defined by $\alpha_f^n(q)$ is the configuration
$p(e'(D(\xi))\cap\{q\}\times\re^n)$, where $p$ is the projection on
$\re^n$, labelled by the map $D(\xi)\buildrel h\over\a D(\gamma^c)\to
MO_c$.

Koschorke and Sanderson [\KSb] use this construction to classify up to
bordism embeddings $e\co M\to Q\times\re^n$ which project to
immersions in $Q$ as the set of homotopy classes $[Q,C_n(MO_c)]$.  To
classify embeddings and knots we need to refine this result to give a
classification up to isotopy.  For this we need to define more general
configuration spaces.  Suppose $X\supset Y$ and the base point $*\in
X$ is not in $Y$.  Define $C_n(X,Y)$ to be based finite subsets of
$\re^n$ with labels in $X$ but with base point (of the configuration)
labelled in $Y$.  The space is topologised in the usual way --- as a
disjoint union of products with equivalence relation defined by
ignoring points labelled by $*$.  Now define a subspace $C_n^0(X)$ of
$C_n(X)$ by restricting configurations to lie in
$\re^n\setminus\{0\}\subset\re^n$.  Then there is a homeomorphism
$C_n(X,Y)\buildrel p_1\times p_2 \over \a C_n^0(X)\times Y$ where
$p_1$ is determined by linear translation of the configuration so that
the base point is translated to $\{0\}$ and $p_2$ is determined by the
base point label of the configuration.

Returning to our immersion $f$ covered by embedding $e$ we have a map
$\beta_f^n\co M\to C_n(MO_c,BO_c)$ given by observing that, for $x\in
M$, the labelled configuration $\alpha_f^nf(x)\in C_n(MO_c)$ has a
distinguished point, given by $e(x)$, which is labelled in $BO_c$.

\rk{Notation}We
shall use the short notation $\I_n^c$ for $C_n(MO_c)$ and $\E_n^c$ for
$C_n(MO_c,BO_c)$.  This notation is intended to remind us that
$\I_n^c$ classifies immersions (up to bordism) and that $\E_n^c$
corresponds to covering embeddings.\endrk

We now have a pull back square:
$$
\lower 20pt\hbox{$\gamma_f^n\co$}\qquad 
\squarediag 
M{\beta_f^n}{\E_n^c}f{\scriptstyle \rm nat}Q{\alpha_f^n}{\I_n^c}
$$
In particular for {\it any} immersion $f\co M\to Q$ there is a pull back square:
$$
\lower 20pt\hbox{$\gamma_f^\infty\co$}\qquad 
\squarediag 
M{\beta_f^\infty}{\E^c_\infty}f{\scriptstyle \rm nat}Q{\alpha_f^\infty}{\I^c_\infty}
$$
since we can choose an embedding in $Q\times \re^\infty$ covering $f$.
 The choice is unique up to isotopy. Note also that a regular homotopy
 of $f$ may be covered by an isotopy in $Q\times \re^\infty$.  A
 simple application of standard obstruction theory shows that in
 codimension 2 a concordance between immersions can be replaced by a
 regular homotopy, hence using proposition \CTProp\ we have the
 following.

\proc{Embedding Theorem}\key{EmbThm} An immersion $f\co M\to Q$ with 
$c=q-m\geq 2$ is regularly homotopic to an immersion which is covered
 by an embedding in $Q\times \re^n$ if and only if there is a pullback
 diagram
$$
\lower 20pt\hbox{$\gamma_F^\infty\co$}\qquad 
\squarediag 
{M\times I}{\beta_F^\infty}{\E^c_\infty}F{\scriptstyle \rm nat}{Q\times I}{\alpha_F^\infty}{\I_\infty^c}
$$
with $f(x)=F(x,0)$, for all $x \in M$, $\alpha_F(q,1)\in \I_n^c$, 
and $\beta_F(q,1)\in \E_n^c$, for all $q \in Q$.\qed
\endproc 

\rk{Remark}The theorem extends to codimension 1, but only up to
concordance.  To be precise, regular homotopy should be replaced by
immersed concordance in both proposition \CTProp\ and the Embedding
Theorem.  This is seen by using the unparametrised multi-compression
theorem in the proof of \CTProp.\endrk

The condition of pullback homotopy in the theorem is difficult to
interpret homotopically, and does not immediately reduce the embedding
problem to obstruction theory.  For a covering embedding to exist the
map of $Q$ to $\I_\infty^c$ must lift to $\I_n^c$.  By [\KSb]
this is the obstruction to finding a solution up to bordism.  To find
a solution up to isotopy the map of $M$ into $C_\infty^0(MO_c)$ must
lift to $C_n^0(MO_c)$ and this defines new obstructions.  The
algebraic topology corresponding to these new obstructions will be
investigated in a future paper.

We now describe some of the geometry involved in these obstructions.
We start by sharpening our description of the Koschorke--Sanderson map
$Q\to \I_n^c$ induced by an immersion $f\co M\to Q$.  Suppose that
this immersion is self-transverse [10].  Then we have $k$-tuple point
manifolds $M_k$ and $Q_k$ and a $k$-fold covering $f_k\co M_k\to Q_k$,
for $k\ge 1$, defined as follows.
$$Q_k=\{X|X\subset M, |X|=k \hbox{ and } |fX|=1\}, 
M_k=\{(X,x)|x\in X \hbox{ and } X\in Q_k\}. $$ 
Further there is a commutative diagram:
$$
\lower 20pt\hbox{$\gamma_f\co$}\qquad 
\upsquarediag 
MfQ{g_k}{h_k}{M_k}{f_k}{Q_k}
$$
where $h_k(X)=f(x)$, for any $x\in X$, $g_k(X,x)=x$, and $f_k(X,x)=X$.

If $f$ is an embedding there is a map $M \to BO_c$ classifying its
normal bundle and a Thom construction $Q\to MO_c$ where $MO_c$ is the
universal Thom space with $c=q-m$.  The Koschorke--Sanderson map $Q\to
\I_\infty^c$ induced by an immersion is a generalisation.  The immersion
$h_k$ has a normal bundle classified by a map $Q_k\to
E\Sigma_k\times_{\Sigma_k}BO_c^{[k]}$ where $[k]$ indicates $k$-fold
Cartesian product and $\Sigma_k$ denotes the symmetric group on $k$
elements.  The immersion $h_k$is an embedding away from $>k$ multiple
points.  There is a corresponding Thom construction. We can replace
$E\Sigma_k$ by $C(k,\re^\infty)$ --- configurations of $k$ distinct
(ordered) points in $\re^\infty$.  A point in the space
$C(k,\re^\infty)\times_{\Sigma_k}BO_c^{[k]}$ may be regarded as a set
of k distinct points in $\re^\infty$ each with a label in $BO_c$.  Now
the configuration space $C_\infty(X)$ (of points of $\re^\infty$
labelled in $X$) is $\coprod_k C(k,\re^\infty)
\times_{\Sigma_k}X^{[k]}/\sim$ where the equivalence relation is given
by ignoring points labelled at the base point.  The immersion $f$ now
determines a map $Q\to\I_\infty^c=C_\infty(MO_c)$, which `puts
together' the Thom constructions for each of the multiple point
manifolds. It is well defined up to homotopy and is determined by the
choice of points in $\re^n$ ie by a commutative diagram:
\commdiag Me{Q\times\re^\infty}f{\rm proj}Q
where $e$ is an embedding.

This map is `transverse' to each of the bases of the various products
of Thom complexes and the pull-backs are the multiple point manifolds.
The process reverses.  Any map $Q\to\I_\infty^c$ can be made transverse in
this way and thus determines a self-transverse immersion.

There is a similar description for the Koschorke--Sanderson map $Q\to
\I_n^c$ determined by an immersion covered by an embedding in 
$Q\times\re^n$.

The multiple point manifolds (with all the normal bundle information)
can be regarded as higher `Hopf invariants' [\KSb].  There are
analogous interpretations for the maps $M\to\E_n^c$ or $\E_\infty^c$
corresponding to the covering embeddings in embeddings in
$Q\times\re^n$ or $Q\times\re^\infty$.  We describe the $\E_n^c$ case;
the $\E_\infty^c$ is similar.

For $1\leq n\leq\infty$ there is a $k$-fold covering $p_k\co
C^\bullet(k,\re^n)\to C(k,\re^n)$ where $C^\bullet(k,\re^n)$ denotes
based sets $A\subset \re^n$, with $|A|=k$.  We can write the space
$C_n^0(X)$ as $\coprod_k C^\bullet(k,\re^n)
\times_{\Sigma_{k-1}}X^{[k-1]}/\sim$ and in case $(X,Y)=(MO_c,BO_c)$
the map $proj\circ\beta^n_f\co M\to C_n(X,Y)=\E_n^c$ puts together the
Thom constructions for the $g_k$.  Using the homeomorphism
$\E_n^c\cong C^0_n(MO_c)\times BO_c$ we can see that this map
classifies (a) the self-transverse system of multiple point sets in
$M$ and (b) the normal bundle on $M$ in $Q$.  Thus part of the second
obstruction (to lifting $M\to \E_\infty^c$ to $\E_n^c$) is given by
higher Hopf invariants for the multiple point sets in $M$.

\sh{Classification of knots}

There is a similar analysis for the knot problem.  Rather than
continuing to consider embeddings in $Q\times\re^n$ which cover
immersions in $Q$, we shall consider embeddings equipped with $n$
independent normal vector fields (for example framed embeddings)
since by the compression theorem these are equivalent.  Such an
embedding $f$ determines a pull-back square:
$$
\lower 20pt\hbox{$\rho_f\co$}\qquad 
\squarediag 
{M}{\beta_{f_1}}{\E_n^c}{f_1}{\scriptstyle \rm nat}{Q}{\alpha_{f_1}}{\I_n^c}
$$
where $f_1={\rm proj}\circ f$.

The arguments used to prove theorem 4.3 now prove:

\proc{Theorem}Suppose that $f,g\,\co M\to Q\times\re^n$
are embeddings equipped with $n$ independent normal vector fields
and that $q-m\ge 2$.  Then $f$ is isotopic to $g$
if and only if $\rho_f$ is homotopic to $\rho_g$ by a pull-back
homotopy. \qed

\proc{Corollary}{\rm (Classification of knots)}\qua Suppose that $q-m\ge2$.
There is a bijection
between isotopy classes of embeddings $f$ of $M$ 
in $Q\times\re^n$ equipped with $n$
independent normal vector fields and pull-back homotopy classes
of squares $\rho_f$. \qed

The corollary gives many knot invariants, for example any homotopy
invariant of $\I_n^c$, pulled back to $Q$, or of $\E_n^c$, pulled back
to $M$.  The former are invariants of the cobordism class of the knot
and the latter are new invariants.  These both contain higher Hopf
invariants (suitably generalised) as outlined above for the embedding
problem.

\rk{Remarks}(1)\qua The last two results also hold in codimension 1, 
but with regular homotopy replaced by regular concordance (see the
remark below theorem \ref{EmbThm}).

(2)\qua The invariants and obstructions discussed above have strong
connections with many existing invariants.  The case $c=1,n=1$ is
studied by Fenn, Rourke and Sanderson see in particular [\FRS;
sections 2 and 4], where classifying spaces related to $\I^1_1$ but
depending of the fundamental rack, are also considered.  There are
combinatorial invariants defined in this case, for example the
generalised James--Hopf invariants, which link with the higher Hopf
invariants described above.

We shall give more details here and also explain connections with
other known obstructions and invariants in a subsequent paper.

\section{Controlling singularities of a projection}

\rk{Definition}A {\sl weakly stratified set} is a set $X$ with
a flag of closed subsets
$$X=S_0\supset S_1\supset S_2\ldots\supset S_t\supset S_{t+1}=\emptyset$$
such that, for each $i=0,\ldots,t$, $S_i-S_{i+1}$ is a manifold.

We also say that $S_i$ is a {\sl weak stratification} of $X$ and we
call the manifolds $S_i-S_{i+1}$ the {\sl strata}. 

\rk{Remarks}This is very much weaker than the usual notion of a 
stratified set --- there is no condition on the neighbourhood of
$S_{i+1}$ in $S_i$ or any relationship between the dimensions of the
strata.

\rk{Definitions}Suppose that $X$ is a weakly stratified set and that 
$X\subset W$ (a manifold).  Suppose that $\xi^n$ is a plane field on
$W$ (ie a $n$--subbundle of $TW$) defined at $X$.  We say $\xi$ is
{\sl weakly normal} to $X$ if $\xi$ is normal to $S_i-S_{i+1}$ for
each $i=0,\ldots,p$.

Suppose that $M^m\subset W^w$ and that $\xi^n$ is a plane field
defined at $M$.  We say that $\xi$ has {\sl regular singularities on
$M$} if it is normal to a weak stratification of $M$.

\proc{Example}\key{VarEx}\rm A plane field in general position has 
regular singularities.  However so do many plane fields which are far
from general position.  Here is an explicit example with $n=1$
constructed by Varley [\Var].  Let $\alpha$ be the $x$--axis in
$\re^3$ and $C\subset\alpha$ a cantor set.  Let $\pi$ be the surface
(a smooth plane) given by $z=y^3$ which contains $\alpha$ and has a
line of inflection along $\alpha$.  Let $\xi$ be the line field
parallel to the $y$--axis and distort $\xi$ to have a small negative
$z$--component off $C$.  Then $\xi$ is tangent to $\pi$ precisely at
$C$ and very far from general position.  But it is has regular
singularities by choosing $\alpha$ as one stratum and $\pi$ as the
next.

It is easy to construct plane fields with non-regular singularities:
for an example with $n=1$, choose any line field for which part of a
flow line lies in $M$.

\proc{$C^0$--Singularity Theorem}
Suppose that $M^m\subset W^w$ and that $\xi$ is an integrable
$n$--plane field on $TW$ defined on a neighbourhood $U$ of $M$ such
that $\xi$ has regular singularities on $M$ and that $n+m<w$.  Suppose
given $\ep >0$ and a homotopy of $\xi$ through integrable plane fields
on $TW$ defined on $U$ finishing with the plane field $\xi'$.  Then
there is an ambient isotopy of $M$ in $W$ which moves points at most
$\ep$ moving the pair $M,\xi|M$ to $M',\xi'|M'$.\endproc

The theorem gives an answer to the $C^0$--Singularity Problem (stated
in section 1) in the case that ${\rm dim}(M)<{\rm dim}(Q)$.  Take
$\xi'$ to be the tangent bundle to the fibres of $p$ then the
singularities of $p|M$ can be made to coincide with those of any plane
field homotopic to $\xi'$ with regular singularities.

Some condition on the singularities is clearly necessary, for example,
again with $n=1$, suppose that part of a flow line of the normal line
field lies in $M$ and is not already vertical (thinking of $\xi'$ as
vertical) then no {\it small} isotopy can make this field vertical.
The condition of regularity is very weak as can be seen by considering
examples similar to \ref{VarEx}.  The condition that the plane fields
are {\it integrable} is needed for our proof, but we do not have an
example to show that it is necessary for the result.  Note that
example \ref{VarEx} is integrable and see [\Var] for more examples.

\prf We shall prove a more general result: $M$ is replaced by
any weakly stratified set $X$ such that dimensions of strata are
$<w-n$ and $\xi$ is weakly normal to $X$.  The idea is to apply the
Normal Deformation Theorem to each stratum in turn starting with $S_t$
and continuing with $S_{t-1}-S_t$ etc).  After the $t-i^{\rm th}$ move
$\xi$ concides with $\xi'$ on $S_i$ and with care this can also be
assumed to be true near $S_i$ and in particular in a neighbourhood of
$S_i$ in $S_{i-1}$.  The next move is made relative to a smaller
neighbourhood, and the result is proved in $t+1$ steps.

So the only point that needs work is the point that $\xi$ can be
assumed to concide with $\xi'$ {\it near} $S_i$.  This is where
integrability is needed.  By integrability we can assume that $\xi$ is
the tangent bundle to a foliation $\cal F$ of $W$ defined near $S_i$
and similarly $\xi'$ is the tangent bundle to $\cal F'$.  Now $\cal F$
and $\cal F'$ coincide at $S_i$ and by a $C^\infty$--small isotopy
$\cal F$ can be moved to $\cal F'$ near $S_i$ and this carries $\xi$
into coincidence with $\xi'$ near $S_i$ as required.\endprf

\rk{Addenda}The proof works (indeed was given) for a weakly stratified
set instead of a submanifold and has a natural relative version
directly from the proof.  There is also a parametrised version which
follows by combining the proof of the parametrised Normal Deformation
Theorem [\CompI, page 425] or [\CompII, page 439] with the last proof:

{\sl Suppose that we have a family of embeddings of $M_t$ in $W$ where
$t\in K$ (a parameter manifold) together with integrable plane fields
$\xi_t$ for $t\in K$ having regular singularities with $M_t$.  Suppose
further that the singularities are ``locally constant'' over $K$ (ie
the whole situation is locally trivial).  Suppose given a
$K$--parameter homotopy of $\xi_t$ through integrable plane fields to
$\xi'_t$.  Then there is a $K$--parameter family of small ambient
isotopies carrying $M_t,\xi_t$ to $M'_t,\xi'_t$.}

\rk{Comments}

(1)\qua Varley [\Var; chapter 3] gives a more general parametrised
theorem in which the singularities are not assumed locally constant.
He assumes that $\xi$ has regular singularities with the parametrised
family as a whole (giving a weak stratification of $M\times K$) and
further that the tangent bundle along the fibres has regular
singularities with respect to the each stratum of this stratification.

(2)\qua The philosophy of our solution to the $C^0$--singularity
problem is that the singularities of an arbitrary plane field are
instrinsic to the embedding $M\subset W$ and invariant under isotopy
of the situation.  Notice that the theorem completely controls the
type of singularity which can arise. 

(3)\qua It is also worth commenting that Eliashberg [\Eliash] has
proved strong results about controlling singularities under more
regular assumptions.  

(4)\qua There is a very nice description in the metastable
range where there is a single singularity manifold (in general
position).  This can be regarded as a fine bordism class in the sense
of Koschorke [\Kos] and is an instrinsic embedded characteristic
class.  For more detail see [\Var; chapter 5].

(5)\qua In future papers (to be written jointly with Varley) we intend
to provide further applications of this theorem to give algebraic
topological information about the existence of singularities of
specified types.
  
%
\references \Addresses
\bye